\documentclass[11pt]{article}
\usepackage{}
\usepackage{mathrsfs}
\usepackage{amsfonts}

\usepackage{CJK}
\usepackage{amsfonts,amssymb,amsmath,mathrsfs}

\usepackage{color,latexsym,amsfonts}
 \newtheorem{theorem}{Theorem}[section]
 \newtheorem{lemma}[theorem]{Lemma}
 \newtheorem{corollary}[theorem]{Corollary}
 \newtheorem{proposition}[theorem]{Proposition}
 
 \newtheorem{Definition}[theorem]{Definition}
 \newtheorem{remark}[theorem]{Remark}
 \newtheorem{condition}[theorem]{Condition}

 \def\blemma{\begin{lemma}\sl{}\def\elemma{\end{lemma}}}
 \def\btheorem{\begin{theorem}\sl{}\def\etheorem{\end{theorem}}}
 \def\bcorollary{\begin{corollary}\sl{}\def\ecorollary{\end{corollary}}}

 \def\bremark{\begin{remark}\sl{}\def\eremark{\end{remark}}}

 \def\beqlb{\begin{eqnarray}}\def\eeqlb{\end{eqnarray}}
 \def\beqnn{\begin{eqnarray*}}\def\eeqnn{\end{eqnarray*}}

 \def\<{\langle}\def\>{\rangle}

 \def\ar{&\!\!}

 \def\eqref#1{{\rm(\ref{#1})}}

\def\d{\textup{d}}

\def\e{\textup{e}}

\def\fin{\hfill$\square$}

\def\newdot{{\kern.8pt\cdot\kern.8pt}}

\def\R{\mathbb{R}}
\def\E{\mathbb{E}}
\def\P{\mathbb{P}}

\def\<{\langle}
\def\>{\rangle}
\def\Proof.{\noindent{\bf Proof.}}

\begin{document}

\

\noindent{}

\bigskip\bigskip

\centerline{\Large\bf Harnack Type Inequalities and Applications for SDE Driven by}

\smallskip

\centerline{\Large\bf Fractional Brownian Motion\footnote{Supported by the Research project of Natural Science Foundation of
   Anhui Provincial Universities (Grant No. KJ2013A134).}
}

\smallskip
\
\bigskip\bigskip

\centerline{Xi-Liang Fan}

\bigskip

\centerline{Department of Mathematics, Anhui Normal University,}

\centerline{Wuhu 241003, China}

\smallskip

\bigskip\bigskip

{\narrower{\narrower

\noindent{\bf Abstract.} For stochastic differential equation driven by  fractional Brownian motion with Hurst parameter $H>1/2$,
 Harnack type inequalities are established by constructing a coupling with unbounded time-dependent drift.
 These inequalities are applied to the study of existence and uniqueness of invariant measure for a discrete Markov semigroup
 constructed in terms of the distribution of the solution.
 Furthermore, we show that entropy-cost inequality holds for the invariant measure.

}}


\bigskip



\section{Introduction}

\setcounter{equation}{0}

The dimensional-free Harnack inequality with powers introduced in \cite{Wang97a}
and the log-Harnack inequality introduced in \cite{Rockner&Wang10a,Wang10a} have been intensively investigated in the context of Markov processes,
see, for example, \cite{Guillin&Wang11a,Liu09a,Ouyang09a,Ouyang&Rockner&Wang12a,Wang07a,Wang11b,Wang&Yuan11a,Zhang13a} and references within.
Harnack type inequalities have become a useful tool in stochastic analysis.
One can see, for instance, \cite{Rockner&Wang03a,Rockner&Wang03b,Wang99a} for strong Feller property and contractivity properties; \cite{Aida&Zhang01a,Aida&Zhang02a,Kawabi05a} for short times behaviors of infinite dimensional diffusions;
\cite{Bobkov&Gentil&Ledoux01a,Gong&Wang01a} for heat kernel estimates, entropy-cost inequalities and transportation cost inequalities.

In this note, we are concerned with stochastic differential equations (SDEs for short) driven by fractional Brownian motion,
whose noise is not Markovian and even more not semimartingale.
Based on the theory of rough path analysis introduced in \cite{Lyons98a},
Coutin and Qian \cite{Coutin&Qian02a} presented an existence and uniqueness result with Hurst parameter $H\in(1/4,1/2)$.
Following the approach of \cite{Zahle98a}, Nualart and R\u{a}\c{s}canu \cite{Nualart&Rascanu02a} derived the existence and uniqueness result with $H>1/2$.
In the previous papers \cite{Fan13a} and \cite{Fan12a,Fan13b}, by using the method of derivative formulae
we have established Harnack type inequalities for SDEs with fractional noises for $H<1/2$ and $H>1/2$, respectively.
Motivated by the work \cite{Wang11b}, where a new technique is applied to construct the coupling for a diffusion process with multiplicative noise,
we will establish directly Harnack type inequalities for SDEs driven by fractional motion with Hurst parameter $H>1/2$.
That is the main purpose of this paper.

This paper is organized as follows.
In the next section, we recall some basic results about fractional integrals and derivatives and fractional Brownian motion.
In section 3, by means of the coupling and Girsanov transformation argument,
the dimension-free Harnack type inequalities and the strong Feller property are shown for SDEs driven by fractional Brownian motion with $H>1/2$.
In terms of the distribution of the solution, we construct a discrete Markov semigroup.
As applications of the inequalities, the existence and uniqueness of invariant probability measure for the corresponding semigroup is proved
and its entropy-cost inequality is established.

\section{Preliminaries}

\setcounter{equation}{0}
\subsection{Fractional integrals and derivatives}

Let $a,b\in\R$ with $a<b$.
For $f\in L^1(a,b)$ and $\alpha>0$, the left- and right-sided fractional Riemann-Liouville integral of $f$ of order $\alpha$
on $[a,b]$ is given by
\beqnn
I_{a+}^\alpha f(x)=\frac{1}{\Gamma(\alpha)}\int_a^x\frac{f(y)}{(x-y)^{1-\alpha}}\d y \eeqnn
and
\beqnn
I_{b-}^\alpha f(x)=\frac{(-1)^{-\alpha}}{\Gamma(\alpha)}\int_x^b\frac{f(y)}{(y-x)^{1-\alpha}}\d y,
\eeqnn
where $x\in(a,b)$ a.e., $(-1)^{-\alpha}=\e^{-i\alpha\pi},\Gamma$ denotes the Euler function.
They extend the usual $n$-order iterated integrals of $f$ for $\alpha=n\in\mathbb{N}$.

Let $I_{a+}^\alpha(L^p)$ (resp. $I_{b-}^\alpha(L^p)$) be the image of $L^p(a,b)$ by the operator $I_{a+}^\alpha$ (resp. $I_{b-}^\alpha$).
If $f\in I_{a+}^\alpha(L^p)$ (resp. $I_{b-}^\alpha(L^p)$) and $0<\alpha<1$,
the function $\phi$ satisfying $f=I_{a+}^\alpha\phi$ (resp. $f=I_{b-}^\alpha\phi$) is unique in $L^p(a,b)$ and it coincides with the left-sided (resp. right-sided) Riemann-Liouville derivative
of $f$ of order $\alpha$ defined by
\beqnn
D_{a+}^\alpha f(x)=\frac{1}{\Gamma(1-\alpha)}\frac{\d}{\d x}\int_a^x\frac{f(y)}{(x-y)^\alpha}\d y\quad
\left(\mbox{resp.}\ D_{b-}^\alpha f(x)=\frac{(-1)^{1+\alpha}}{\Gamma(1-\alpha)}\frac{\d}{\d x}\int_x^b\frac{f(y)}{(y-x)^\alpha}\d y\right).
\eeqnn
The corresponding Weyl representation reads as follow
\beqnn
D_{a+}^\alpha f(x)=\frac{1}{\Gamma(1-\alpha)}\left(\frac{f(x)}{(x-a)^\alpha}+\alpha\int_a^x\frac{f(x)-f(y)}{(x-y)^{\alpha+1}}\d y\right)\\
\left(\mbox{resp.}\ \ D_{b-}^\alpha f(x)=\frac{(-1)^\alpha}{\Gamma(1-\alpha)}\left(\frac{f(x)}{(b-x)^\alpha}+\alpha\int_x^b\frac{f(x)-f(y)}{(y-x)^{\alpha+1}}\d y\right)\right),
\eeqnn
where the convergence of the integrals at the singularity $y=x$ holds pointwise for almost all $x$ if $p=1$ and in $L^p$-sense if $1<p<\infty$.

By definition, we have the following inversion formulas \begin{itemize}
\item[]$I_{a+}^\alpha(D_{a+}^\alpha f)=f,\quad \forall f\in  I_{a+}^\alpha(L^p); \quad I_{b-}^\alpha(D_{b-}^\alpha f)=f,\quad \forall f\in  I_{b-}^\alpha(L^p)$;
\item[]$D_{a+}^\alpha(I_{a+}^\alpha f)=f,\quad D_{b-}^\alpha(I_{b-}^\alpha f)=f,\quad \forall f\in L^1(a,b)$.
\end{itemize}

For any $\lambda\in(0,1)$, we denote by $C^\lambda(a,b)$ the set of $\lambda$-H\"{o}lder continuous functions on $[a,b]$.
Recall from \cite{Samko&Kilbas&Marichev} that we have
\begin{itemize}
\item[(i)] if $\alpha<1/p$ and $q=p/(1-\alpha p)$, then $I_{a+}^\alpha(L^p)=I_{b-}^\alpha(L^p)\subset L^q(a,b)$;
\item[(ii)] if $\alpha>1/p$, then $I_{a+}^\alpha(L^p)\cup I_{b-}^\alpha(L^p)\subset C^{\alpha-1/p}(a,b)$.
\end{itemize}

Suppose that $f\in C^\lambda(a,b)$ and $g\in C^\mu(a,b)$ with $\lambda+\mu>1$.
By \cite{Young36a}, the Riemann-Stieltjes integral $\int_a^bf\d g$ exists.
In \cite{Zahle98a}, Z\"{a}hle provides an explicit expression for the integral $\int_a^bf\d g$ in terms of fractional derivative. Let $\lambda>\alpha$ and $\mu>1-\alpha$.
Then the Riemann-Stieltjes integral can be expressed as
\beqlb\label{2.0}
\int_a^bf\d g=(-1)^\alpha\int_a^bD_{a+}^\alpha f(t)D_{b-}^{1-\alpha}g_{b-}(t)\d t,
\eeqlb
where $g_{b-}(t)=g(t)-g(b)$. \\
The relation \eqref{2.0} can be regarded as fractional integration by parts formula.

\subsection{Fractional Brownian motion}

In this subsection, we will recall some results about fractional Brownian motion.
The main references for all these results are \cite{Alos&Mazet&Nualart01a}, \cite{Biagini&Hu08a}, \cite{Decreusefond&Ustunel98a} and \cite{Nualart06a}.

Fix a time interval $[0,T]$.
Let $H\in(0,1)$.
The $d$-dimensional fractional Brownian motion with Hurst parameter $H$ on the probability space $(\Omega,\mathscr{F},\mathbb{P})$
can be defined as the centered Gauss process $B^H=\{B_t^H, t\in[0,T]\}$ with covariance function
$\E B_t^{H,i}B_s^{H,j}=R_H(t,s)\delta_{i,j}$,
where
\beqnn
R_H(t,s)=\frac{1}{2}\left(t^{2H}+s^{2H}-|t-s|^{2H}\right).
\eeqnn
In particular, if $H=1/2, B^H$ is a $d$-dimensional Brownian motion.
Besides, one can show that $\E|B_t^{H,i}-B_s^{H,i}|^p=C(p)|t-s|^{pH},\ \forall p\geq 1$.
Consequently, $B^{H,i}$ have $(H-\epsilon)$-H\"{o}lder continuous paths for all $\epsilon>0,\ i=1,\cdot\cdot\cdot,d$.

For each $t\in[0,T]$, we denote by $\mathcal {F}_t$ the $\sigma$-algebra generated by the random variables $\{B_s^H:s\in[0,t]\}$ and the
$\mathbb{P}$-null sets.

We denote by $\mathscr{E}$ the set of step functions defined on $[0,T]$.
Let $\mathcal {H}$ be the Hilbert space defined as the closure of $\mathscr{E}$ with respect to the scalar product
\beqnn
\langle (I_{[0,t_1]},\cdot\cdot\cdot,I_{[0,t_d]}),(I_{[0,s_1]},\cdot\cdot\cdot,I_{[0,s_d]})\rangle_\mathcal {H}=\sum\limits_{i=1}^dR_H(t_i,s_i).
\eeqnn
The mapping $(I_{[0,t_1]},\cdot\cdot\cdot,I_{[0,t_d]})\mapsto\sum_{i=1}^dB_{t_i}^{H,i}$ can be extended to an isometry
between $\mathcal {H}$ and the Gauss space $\mathcal {H}_1$ associated with $B^H$.
We denote the isometry between $\mathcal {H}$ and $\mathcal {H}_1$ by $\phi\mapsto B^H(\phi)$.
On the other hand, the covariance kernel $R_H(t,s)$ can be written as
\beqnn
R_H(t,s)=\int_0^{t\wedge s}K_H(t,r)K_H(s,r)\d r,
\eeqnn
where $K_H$ is a square integrable kernel given by
\beqnn
K_H(t,s)=\Gamma\left(H+\frac{1}{2}\right)^{-1}(t-s)^{H-\frac{1}{2}}F\left(H-\frac{1}{2},\frac{1}{2}-H,H+\frac{1}{2},1-\frac{t}{s}\right),
\eeqnn
in which $F(\cdot,\cdot,\cdot,\cdot)$ is the Gauss hypergeometric function (for details see \cite{Nikiforov&Uvarov88}).

Now, we define the linear operator $K_H^*:\mathscr{E}\rightarrow L^2([0,T],\R^d)$ by
\beqnn
(K_H^*\phi)(s)=K_H(T,s)\phi(s)+\int_s^T(\phi(r)-\phi(s))\frac{\partial K_H}{\partial r}(r,s)\d r.
\eeqnn
It can be shown (see \cite{Alos&Mazet&Nualart01a}) that, for all $\phi,\psi\in\mathscr{E}$,
\beqnn
\langle K_H^*\phi,K_H^*\psi\rangle_{L^2([0,T],\R^d)}=\langle\phi,\psi\rangle_\mathcal {H},
\eeqnn
and therefore $K_H^*$ is  an isometry between $\mathcal{H}$ and $L^2([0,T],\R^d)$.
Consequently, the fractional Brownian motion $B^H$ has the following integral representation
\beqnn
 B^H_t=\int_0^tK_H(t,s)\d W_s,
\eeqnn
where $\{W_t=B^H((K_H^*)^{-1}{\rm I}_{[0,t]}),t\in[0,T]\}$ is a Wiener process.

Consider the operator $K_H: L^2([0,T],\R^d)\rightarrow I_{0+}^{H+1/2 }(L^2([0,T],\R^d))$ associated with the integrable kernel $K_H(\cdot,\cdot)$
\beqnn
(K_Hf^i)(t)=\int_0^tK_H(t,s)f^i(s)\d s,\ i=1,\cdot\cdot\cdot,d.
\eeqnn
It can be proved (see \cite{Decreusefond&Ustunel98a}) that $K_H$ is an isomorphism and moreover,  for each $f\in L^2([0,T],\mathbb{R}^d)$,
\beqnn
(K_H f)(s)=I_{0+}^{2H}s^{1/2-H}I_{0+}^{1/2-H}s^{H-1/2}f,\ H\leq1/2,\\
(K_H f)(s)=I_{0+}^{1}s^{H-1/2}I_{0+}^{H-1/2}s^{1/2-H}f,\ H\geq1/2.
\eeqnn
Consequently, for each $h\in I_{0+}^{H+1/2}(L^2([0,T],\R^d))$, the inverse operator $K_H^{-1}$ is of the form
\beqlb\label{2.0'}
(K_H^{-1}h)(s)=s^{H-1/2}D_{0+}^{H-1/2}s^{1/2-H}h',\ H>1/2,
\eeqlb
\beqlb\label{2.0''}
(K_H^{-1}h)(s)=s^{1/2-H}D_{0+}^{1/2-H}s^{H-1/2}D_{0+}^{2H}h,\ H<1/2.
\eeqlb
In particular, if $h$ is absolutely continuous, we can write for $H<1/2$
\beqlb\label{2.0'''}
(K_H^{-1}h)(s)=s^{H-1/2}I_{0+}^{1/2-H}s^{1/2-H}h' .
\eeqlb

In the paper, we are interested in the equation driven by fractional Brownian motion with Hurst parameter $H>1/2$
\beqlb\label{2.1}
\d X_t=b(t,X_t)\d t+\sigma(t)\d B^H_t,\ X_0=x\in\R^d,\ t\in[0,T],
\eeqlb
where $b:[0,T]\times\R^d\rightarrow\R^d,\ \sigma:[0,T]\rightarrow\R^d$.

Our aim is to establish Harnack type inequalities for \eqref{2.1}, and moreover present some applications.
We will also need some notations.
Define $P_tf(x)=\E f(X_t^x), \ t\in[0,T], \ f\in \mathscr{B}_b(\R^d)$,
where $X_t^x$ is the solution of \eqref{2.1} with the initial point $x$ and $\mathscr{B}_b(\R^d)$ is the set of all bounded measurable functions on $\R^d$.
For all $f\in C^\lambda(0,T)$, let $\|f\|_\infty=\sup_{0\leq t\leq T}|f(t)|$ and $\|f\|_\lambda=\sup_{0\leq s<t\leq T}\frac{|f(t)-f(s)|}{|t-s|^\lambda}$.

\section{Harnack type inequalities and their applications}

\setcounter{equation}{0}

We begin with the assumption (H1)
\begin{itemize}
\item[(i)] $b$ is Lipschitz continuous with non-negative constant $K$:
$$|b(t,x)-b(t,y)|\leq K|x-y|,\ \forall t\in[0,T],\ x,y\in\R^d,$$
 and $b(\cdot,x)$ is Lipschitz continuous;
\par
\item[(ii)] $\sigma^{-1}$ is H\"{o}lder continuous of order $H-1/2<\alpha_0\leq1$ with non-negative constant $\bar{K}$:
 $$|\sigma^{-1}(t)-\sigma^{-1}(s)|\leq \bar{K}|t-s|^{\alpha_0},\ \forall t,s\in[0,T],$$
  and $\sigma$ is bounded.
\end{itemize}
It has been shown in \cite{Nualart&Rascanu02a} that under the above assumption, there exists a unique adapted solution to equation \eqref{2.1}
whose trajectories are H\"{o}lder continuous of order $H-\epsilon$ for any $\epsilon>0$.

For this kind of equation, main result reads as follow.

\btheorem\label{T3.1}
Assume (H1).
Then there exist positive constants $C,C'$ and $C''$ such that
\begin{itemize}
\item[(1)] the log-Harnack inequality
\beqnn
P_T\log f(y)\leq\log P_T f(x)+(C+C'T+C''T^{2\alpha_0})\frac{T^{2(1-H)}}{(1-\e^{-\frac{2}{3}K T})^3}|x-y|^2,\ x,y\in\R^d
\eeqnn
holds for any positive $f\in\mathscr{B}_b(\R^d)$;
\item[(2)] the Harnack inequality
\beqnn
(P_T f)^p(y)\leq P_T f^p(x)\exp\left[\frac{p}{p-1}(C+C'T+C''T^{2\alpha_0})\frac{T^{2(1-H)}}{(1-\e^{-\frac{2}{3}K T})^3}|x-y|^2\right],\ x,y\in\R^d
\eeqnn
holds for all non-negative $f\in\mathscr{B}_b(\R^d)$.
\end{itemize}
\etheorem

To prove the theorem, we first construct a coupling equation.

Let $x,y\in R^d$ such that $x\neq y$.
For $\theta_0\in(0,2)$, let
\beqlb\label{3.0}
\zeta(t)=\frac{2-\theta_0}{2K}(1-\e^{\frac{2}{3}K(t-T)}),\ t\in[0,T].
\eeqlb
Then $\zeta$ is smooth, nonincreasing and strictly positive on $[0,T)$ satisfying
\beqlb\label{3.0'}
3\zeta'(t)-2K\zeta(t)+2=\theta_0,\ t\in[0,T].
\eeqlb
Let $X_t$ solve the equation \eqref{2.1} and introduce the coupling equation as follows
\beqlb\label{3.1}
\d Y_t=b(t,Y_t)\d t+\sigma(t)\d B^H_t+\frac{X_t-Y_t}{\zeta(t)}\d t,\ \ Y_0=y\in\R^d.
\eeqlb
By \cite[Theorem 2.1]{Nualart&Rascanu02a},  \eqref{2.1} and \eqref{3.1} have a unique solution $(X_t,Y_t)$ for $t\in[0,T)$.
That is, the fact that the additional drift in  \eqref{3.1} is singular at time $T$ leads to $Y_t$ is well solved only before time $T$.
To solve $Y_t$ for all $t\in[0,T]$, we need to reformulate the equation by using a new fractional Brownian motion.
To this end, let
\beqnn
\tilde{B}^H_t&=&B^H_t+\int_0^t\sigma^{-1}(s)\frac{X_s-Y_s}{\zeta(s)}\d s\cr
&=&\int_0^tK_H(t,s)\left(\d W_s+K_H^{-1}\left(\int_0^\cdot\sigma^{-1}(\theta)\frac{X_\theta-Y_\theta}{\zeta(\theta)}\d\theta\right)(s)\d s\right)\cr
&=:&\int_0^tK_H(t,s)\d \tilde{W}_s, t\in[0,T).
\eeqnn
Now, for $t\in[0,T)$ we set
\beqnn
R(t)&=&\exp\Bigg[-\int_0^t\left\langle K_H^{-1}\left(\int_0^\cdot\sigma^{-1}(\theta)\frac{X_\theta-Y_\theta}{\zeta(\theta)}\d\theta\right)(r),\d W_r\right\rangle\cr
&&-\frac{1}{2}\int_0^t\left| K_H^{-1}\left(\int_0^\cdot\sigma^{-1}(\theta)\frac{X_\theta-Y_\theta}{\zeta(\theta)}\d\theta\right)(r)\right|^2\d r\Bigg].
\eeqnn

\blemma\label{L3.1}
Assume (H1).
Then, there exist positive constants $C,C'$ and $C''$ such that
\beqnn
\E(R(t)\log R(t))\leq(C+C'T+C''T^{2\alpha_0})\frac{T^{2(1-H)}}{(1-\e^{-\frac{2}{3}K T})^3}|x-y|^2,\ \ \forall t\in[0,T).
\eeqnn
Moreover, $R(T):=\lim_{t\uparrow T}R(t)$ exists and $\{R(t)\}_{t\in[0,T]}$ is a uniformly integrable martingale.
\elemma

\emph{Proof.}
Fixed $s_0\in[0,T)$.
Note that
\beqnn
\d(X_t-Y_t)=(b(t,X_t)-b(t,Y_t))\d t-\frac{X_t-Y_t}{\zeta(t)}\d t,
\eeqnn
then we obtain
\beqnn
\d|X_t-Y_t|^2&=&2\langle X_t-Y_t,b(t,X_t)-b(t,Y_t)\rangle\d t-2\frac{|X_t-Y_t|^2}{\zeta(t)}\d t\cr
&\leq&2\left(K-\frac{1}{\zeta(t)}\right)|X_t-Y_t|^2\d t,\ t\leq s_0.
\eeqnn
This, together with \eqref{3.0'}, leads to
\beqnn
\d \frac{|X_t-Y_t|^2}{\zeta^3(t)}&=&\frac{\d |X_t-Y_t|^2}{\zeta^3(t)}-\frac{3\zeta'(t)}{\zeta^4(t)}|X_t-Y_t|^2\d t\cr
&\leq&-\frac{3\zeta'(t)-2K\zeta(t)+2}{\zeta^4(t)}|X_t-Y_t|^2\d t\cr
&=&-\frac{\theta_0}{\zeta^4(t)}|X_t-Y_t|^2\d t,\ t\leq s_0.
\eeqnn
As a consequence, multiplying by $1/\theta_0$ and integrating from $0$ to $s_0$ yield
\beqlb\label{3.3}
\int_0^{s_0}\frac{|X_t-Y_t|^2}{\zeta^4(t)}\d t+\frac{|X_{s_0}-Y_{s_0}|^2}{\theta_0\zeta^3(s_0)}\leq\frac{|x-y|^2}{\theta_0\zeta^3(0)},\ s_0\in[0,T).
\eeqlb
On the other hand, it follows from \eqref{2.0'} that
\beqlb\label{3.4}
&&\int_0^{s_0}\left| K_H^{-1}\left(\int_0^\cdot\sigma^{-1}(\theta)\frac{X_\theta-Y_\theta}{\zeta(\theta)}\d\theta\right)(r)\right|^2\d r\cr
&=&\frac{1}{\Gamma(3/2-H)^2}\int_0^{s_0}\Bigg|r^{\frac{1}{2}-H}\sigma^{-1}(r)\frac{X_r-Y_r}{\zeta(r)}\cr
&&+\left(H-\frac{1}{2}\right)r^{H-\frac{1}{2}}\int_0^{r}\frac{r^{\frac{1}{2}-H}-\theta^{\frac{1}{2}-H}}{(r-\theta)^{\frac{1}{2}+H}}\sigma^{-1}(\theta)\frac{X_\theta-Y_\theta}{\zeta(\theta)}\d\theta\cr
&&+\left(H-\frac{1}{2}\right)\int_0^{r}
\frac{\sigma^{-1}(r)\frac{X_r-Y_r}{\zeta(r)}-\sigma^{-1}(\theta)\frac{X_\theta-Y_\theta}{\zeta(\theta)}}{(r-\theta)^{\frac{1}{2}+H}}\d\theta\Bigg|^2\d r.
\eeqlb
Next, we are to estimate \eqref{3.4}.\\
Note that, by (H1) and \eqref{3.3}, we conclude that
\beqlb\label{3.5}
\int_0^{s_0}\left|r^{\frac{1}{2}-H}\sigma^{-1}(r)\frac{X_r-Y_r}{\zeta(r)}\right|^2\d r
&=&\int_0^{s_0}r^{1-2H}|\sigma^{-1}(r)|^2\zeta(r)\frac{|X_r-Y_r|^2}{\zeta^3(r)}\d r\cr
&\leq&\frac{\|\sigma^{-1}\|_\infty^2\|\zeta\|_\infty}{2(1-H)\zeta^3(0)}T^{2(1-H)}|x-y|^2
\eeqlb
and
\beqlb\label{3.6}
&&\int_0^{s_0}\left|r^{H-\frac{1}{2}}\int_0^{r}\frac{r^{\frac{1}{2}-H}-\theta^{\frac{1}{2}-H}}{(r-\theta)^{\frac{1}{2}+H}}
\sigma^{-1}(\theta)\frac{X_\theta-Y_\theta}{\zeta(\theta)}\d\theta\right|^2\d r\cr
&\leq&\|\sigma^{-1}\|_\infty^2\sup\limits_{\theta\in[0,s_0]}\frac{|X_\theta-Y_\theta|^2}{\zeta^2(\theta)}
\int_0^{s_0}r^{2H-1}\left(\int_0^r\frac{r^{\frac{1}{2}-H}-\theta^{\frac{1}{2}-H}}{(r-\theta)^{\frac{1}{2}+H}}\d\theta\right)^2\d r\cr
&\leq&\frac{(C_0\|\sigma^{-1}\|_\infty)^2\|\zeta\|_\infty}{2(1-H)\zeta^3(0)}T^{2(1-H)}|x-y|^2,
\eeqlb
where we use the relation
$$\int_0^r\frac{r^{\frac{1}{2}-H}-\theta^{\frac{1}{2}-H}}{(r-\theta)^{\frac{1}{2}+H}}\d\theta
=\int_0^1\frac{s^{\frac{1}{2}-H}-1}{(1-s)^{\frac{1}{2}+H}}\d s\cdot r^{1-2H}=:C_0r^{1-2H}.$$
Besides, by \eqref{2.1} and \eqref{3.1}, we have
\beqnn
&&\int_0^{r}
\frac{\sigma^{-1}(r)\frac{X_r-Y_r}{\zeta(r)}-\sigma^{-1}(\theta)\frac{X_\theta-Y_\theta}{\zeta(\theta)}}{(r-\theta)^{\frac{1}{2}+H}}\d\theta\cr
&=&\frac{\sigma^{-1}(r)(X_r-Y_r)}{\zeta(r)}\int_0^{r}
\frac{1}{(r-\theta)^{\frac{1}{2}+H}}\frac{\zeta(\theta)-\zeta(r)}{\zeta(\theta)}\d\theta\cr
&&+\int_0^{r}
\frac{1}{(r-\theta)^{\frac{1}{2}+H}}\frac{\sigma^{-1}(r)(X_r-Y_r)-\sigma^{-1}(\theta)(X_\theta-Y_\theta)}{\zeta(\theta)}\d\theta.\cr
&=&\frac{\sigma^{-1}(r)(X_r-Y_r)}{\zeta(r)}\int_0^{r}
\frac{1}{(r-\theta)^{\frac{1}{2}+H}}\frac{\zeta(\theta)-\zeta(r)}{\zeta(\theta)}\d\theta\cr
&&+\int_0^{r}\frac{\sigma^{-1}(r)-\sigma^{-1}(\theta)}{(r-\theta)^{\frac{1}{2}+H}}\frac{X_\theta-Y_\theta}{\zeta(\theta)}\d\theta\cr
&&+\sigma^{-1}(r)\int_0^{r}\frac{1}{(r-\theta)^{\frac{1}{2}+H}\zeta(\theta)}\int_\theta^{r}[b(s,X_s)-b(s,Y_s)]\d s\d\theta\cr
&&-\sigma^{-1}(r)\int_0^{r}\frac{1}{(r-\theta)^{\frac{1}{2}+H}\zeta(\theta)}\int_\theta^{r}\frac{X_s-Y_s}{\zeta(s)}\d s\d\theta\cr
&=:&J_1(r)+J_2(r)+J_3(r)+J_4(r).
\eeqnn
Observe that, by the definition of $\zeta$ and \eqref{3.3}, we get
\beqlb\label{3.7}
\int_0^{s_0}|J_1(r)|^2\d r
&\leq&\|\sigma^{-1}\|_\infty^2\int_0^{s_0}\frac{|X_r-Y_r|^2}{\zeta^4(r)}\left(\int_0^r\frac{\zeta(\theta)-\zeta(r)}{(r-\theta)^{\frac{1}{2}+H}}\d\theta\right)^2\d r\cr
&=&\left[\frac{(2-\theta_0)\|\sigma^{-1}\|_\infty
\e^{-\frac{2K}{3}T}}{2K}\right]^2\int_0^{s_0}\frac{|X_r-Y_r|^2}{\zeta^4(r)}
  \left(\int_0^r\frac{\e^{\frac{2K}{3} r}-\e^{\frac{2K}{3}\theta}}{(r-\theta)^{\frac{1}{2}+H}}\d\theta\right)^2\d r\cr
&=&\left[\frac{(2-\theta_0)\|\sigma^{-1}\|_\infty
\e^{-\frac{2K}{3} T}}{3}\right]^2\int_0^{s_0}\frac{|X_r-Y_r|^2}{\zeta^4(r)}
  \left(\int_0^re^{\frac{2K}{3}s}\d s\int_0^s\frac{1}{(r-\theta)^{\frac{1}{2}+H}}\d\theta\right)^2\d r\cr
&\leq&\left[\frac{(2-\theta_0)\|\sigma^{-1}\|_\infty}{3(H-1/2)(3/2-H)}\right]^2\int_0^{s_0}\frac{|X_r-Y_r|^2}{\zeta^4(r)}r^{3-2H}\d r\cr
&\leq&\left[\frac{(2-\theta_0)\|\sigma^{-1}\|_\infty}{3(H-1/2)(3/2-H)}\right]^2\frac{1}{\theta_0\zeta^3(0)}T^{3-2H}|x-y|^2.
\eeqlb
and with the help of (H1) and \eqref{3.3}, it follows \beqlb\label{3.8}
\int_0^{s_0}|J_2(r)|^2\d r
&\leq&\sup_{0\leq\theta\leq s_0}\frac{|X_\theta-Y_\theta|^2}{\zeta^2(\theta)}\int_0^{s_0}\left|\int_0^r\frac{\sigma^{-1}(r)-\sigma^{-1}(\theta)}{(r-\theta)^{\frac{1}{2}+H}}\d\theta\right|^2\d r\cr
&\leq&\frac{\bar{K}^2\|\zeta\|_\infty}{2(\alpha_0-H+1)(\alpha_0-H+1/2)^2\zeta^3(0)}T^{2(\alpha_0-H+1)}|x-y|^2.
\eeqlb
In view of (H1), the Fubini theorem, the Cauchy-Schwarz inequality and \eqref{3.3}, we get
\beqlb\label{3.9}
\int_0^{s_0}|J_3(r)|^2\d r
&\leq&(K\|\sigma^{-1}\|_\infty)^2\int_0^{s_0}\left(\int_0^r|X_s-Y_s|\d s\int_0^s\frac{1}{\zeta(\theta)(r-\theta)^{\frac{1}{2}+H}}\d\theta\right)^2\d r\cr
&\leq&\left(\frac{K\|\sigma^{-1}\|_\infty}{H-1/2}\right)^2
      \int_0^{s_0}\left(\int_0^r\frac{|X_s-Y_s|}{\zeta(s)}(r-s)^{\frac{1}{2}-H}\d s\right)^2\d r\cr
&\leq&\left(\frac{K\|\sigma^{-1}\|_\infty\|\zeta\|_\infty}{H-1/2}\right)^2
      \int_0^{s_0}\int_0^r\frac{|X_s-Y_s|^2}{\zeta^4(s)}\d s\int_0^r(r-s)^{1-2H}\d s\d r\cr
&\leq&\left(\frac{K\|\sigma^{-1}\|_\infty\|\zeta\|_\infty}{H-1/2}\right)^2
      \frac{1}{2(1-H)(3-2H)\theta_0\zeta^3(0)}T^{3-2H}|x-y|^2.
\eeqlb
As for $J_4(r)$, similar to $J_3(r)$ we have
\beqlb\label{3.10}
\int_0^{s_0}|J_4(r)|^2\d r
\leq\left(\frac{\|\sigma^{-1}\|_\infty}{H-1/2}\right)^2\frac{1}{2(1-H)(3-2H)\theta_0\zeta^3(0)}T^{3-2H}|x-y|^2.
\eeqlb
Substituting \eqref{3.5}, \eqref{3.6}, \eqref{3.7}, \eqref{3.8}, \eqref{3.9}, \eqref{3.10} into \eqref{3.4}, we conclude that
\beqlb\label{3.11}
&&\frac{1}{2}\int_0^{s_0}\left| K_H^{-1}\left(\int_0^\cdot\sigma^{-1}(\theta)\frac{X_\theta-Y_\theta}{\zeta(\theta)}\d\theta\right)(r)\right|^2\d r\cr
&\leq&\frac{3}{\Gamma(3/2-H)^2\zeta^3(0)}
\Bigg\{\frac{\|\sigma^{-1}\|_\infty^2\|\zeta\|_\infty}{2(1-H)}
+\frac{[C_0(H-1/2)\|\sigma^{-1}\|_\infty]^2\|\zeta\|_\infty}{2(1-H)}\cr
&&+\left[\frac{(2-\theta_0)\|\sigma^{-1}\|_\infty}{3(3/2-H)}\right]^2\frac{1}{\theta_0}T
+\frac{[\bar{K}(H-1/2)]^2\|\zeta\|_\infty}{2(\alpha_0-H+1)(\alpha_0-H+1/2)^2}T^{2\alpha_0}\cr
&&+\frac{\left(K\|\sigma^{-1}\|_\infty\|\zeta\|_\infty\right)^2}{2(1-H)(3-2H)\theta_0}T
+\frac{\|\sigma^{-1}\|_\infty^2}{2(1-H)(3-2H)\theta_0}T\Bigg\}T^{2(1-H)}|x-y|^2\cr
&=:&(C+C'T+C''T^{2\alpha_0})\frac{T^{2(1-H)}}{(1-\e^{-\frac{2}{3}K T})^3}|x-y|^2.
\eeqlb
Then it follows that
\beqnn
\E\exp\left[\frac{1}{2}\int_0^{s_0}\left| K_H^{-1}\left(\int_0^\cdot\sigma^{-1}(\theta)\frac{X_\theta-Y_\theta}{\zeta(\theta)}\d\theta\right)(r)\right|^2\d r\right]<\infty.
\eeqnn
Consequently, $\{R_t\}_{t\in[0,s_0]}$ is a martingale and $\{\tilde{W}_t\}_{t\in[0,s_0]}$ is a $d$-dimensional Brownian motion under $R(s_0)\d\P$.
Note that by the definition of $\tilde{W}$, we have
\beqnn
\log R(t)&=&-\int_0^t\left\langle K_H^{-1}\left(\int_0^\cdot\sigma^{-1}(\theta)\frac{X_\theta-Y_\theta}{\zeta(\theta)}\d\theta\right)(r),\d \tilde{W}_r\right\rangle\cr
&&+\frac{1}{2}\int_0^t\left| K_H^{-1}\left(\int_0^\cdot\sigma^{-1}(\theta)\frac{X_\theta-Y_\theta}{\zeta(\theta)}\d\theta\right)(r)\right|^2\d r.
\eeqnn
Combining this with \eqref{3.11} we obtain
\beqnn
\E(R(s_0)\log R(s_0))=\E_{s_0}\log R(s_0)\leq(C+C'T+C''T^{2\alpha_0})\frac{T^{2(1-H)}}{(1-\e^{-\frac{2}{3}K T})^3}|x-y|^2,\ \ s_0\in[0,T),
\eeqnn
where $\E_{s_0}$ denotes the expectation under the probability $R(s_0)\d\P$.\\
Hence, $\{R_t\}_{t\in[0,T)}$ is a uniformly integrable martingale.
As a consequence, by the martingale convergence theorem, we know that
$R(T):=\lim_{t\uparrow T}R(t)$ exists and $\{R(t)\}_{t\in[0,T]}$ is a uniformly integrable martingale.
This completes the proof.
\fin

Lemma \ref{L3.1} ensures that $\{\tilde{B}_t\}_{t\in[0,T]}$ is a $d$-dimensional fractional Brownian motion under the
probability $R(T)\d\P$ by the Girsanov theorem for the fractional Brownian motion (see e.g. \cite[Theorem 4.9]{Decreusefond&Ustunel98a} or \cite[Theorem 2]{Nualart&Ouknine02b}), and together with the Fatou lemma, there holds
\beqlb\label{3.12}
\E(R(T)\log R(T))\leq(C+C'T+C''T^{2\alpha_0})\frac{T^{2(1-H)}}{(1-\e^{-\frac{2}{3}K T})^3}|x-y|^2.
\eeqlb
Rewrite \eqref{2.1} and \eqref{3.1} as follow
\beqlb\label{3.13}
\d X_t=b(t,X_t)\d t+\sigma(t)\d\tilde{B}^H_t-\frac{X_t-Y_t}{\zeta(t)}\d t,\ \ X_0=x,
\eeqlb
\beqlb\label{3.14}
\d Y_t=b(t,Y_t)\d t+\sigma(t)\d\tilde{B}^H_t,\ \ Y_0=y.
\eeqlb
As a consequence, $Y_t$ can be solved from the equation \eqref{3.14} up to time $T$.
Note that the relation $\int_0^T1/\zeta(t)\d t=\infty$ holds, we shall see that the
coupling $(X_t,Y_t)$ is successful up to time $T$.
Thus, $X_T=Y_T$ holds $R(T)\d\P$ a.s. and for the same initial points the distribution of $Y_T$ under $R(T)\d\P$ coincides with
that of $X_T$ under $\P$.
Therefore, we will derive the desired Harnack type inequalities.

\emph{Proof of Theorem \ref{T3.1}.}
Letting $s_0\uparrow T$ in \eqref{3.3}, we conclude that
\beqlb\label{3.15}
\int_0^T\frac{|X_t-Y_t|^2}{\zeta^4(t)}\d t\leq\frac{|x-y|^2}{\theta_0\zeta^3(0)}.
\eeqlb
This implies that the coupling time $\tau:=\inf\{t\in[0,T]:X_t=Y_t\}\leq T$ and so, $X_T=Y_T$ holds, where we set $\inf\emptyset=\infty$ by convention.
Indeed, if there exists $\omega\in\Omega$ such that $\tau(\omega)>T$, then the continuity of the processes $X$ and $Y$ yields
\beqnn
\inf\limits_{t\in[0,T]}|X_t(\omega)-Y_t(\omega)|>0.
\eeqnn
As a consequence, we get
\beqnn
\int_0^T\frac{|X_t(\omega)-Y_t(\omega)|^2}{\zeta^4(t)}\d t
&\geq&\inf\limits_{t\in[0,T]}|X_t(\omega)-Y_t(\omega)|^2\int_0^T\frac{1}{\zeta^4(t)}\d t\cr
&\geq&\inf\limits_{t\in[0,T]}|X_t(\omega)-Y_t(\omega)|^2\frac{1}{T^3}\left(\int_0^T\frac{1}{\zeta(t)}\d t\right)^4=\infty.
\eeqnn
This contradicts with \eqref{3.15}.\\
Therefore, by the Young inequality and the H\"{o}lder inequality we obtain
\beqlb\label{3.16}
P_T\log f(y)=\E(R(T)\log f(Y_T^y))=\E(R(T)\log f(X_T^x))\leq\E(R(T)\log R(T))+\log P_T f(x),
\eeqlb
\beqlb\label{3.17}
(P_T f)^p(y)=\E(R(T)f(Y_T^y))^p=\E(R(T)f(X_T^x))^p\leq P_T f^p(x)\left(\E R(T)^{\frac{p}{p-1}}\right)^{p-1},
\eeqlb
where the superscripts $x$ and $y$ stand for the initial points of corresponding equations, respectively.\\
Combining \eqref{3.16} and \eqref{3.12} implies the desired log-Harnack inequality.
As for the Harnack inequality, by the definition of $R(t)$ and \eqref{3.11}, we have, for $t\in[0,T)$,
\beqnn
&&\E R(t)^{\frac{p}{p-1}}=\E_t R(t)^{\frac{1}{p-1}}\cr
&=&\E_t\exp\Bigg[-\frac{1}{p-1}\int_0^t\left\langle K_H^{-1}\left(\int_0^\cdot\sigma^{-1}(\theta)\frac{X_\theta-Y_\theta}{\zeta(\theta)}\d\theta\right)(r),\d\tilde{W}_r\right\rangle\cr
&&+\frac{1}{2}\frac{1}{p-1}\int_0^t\left| K_H^{-1}\left(\int_0^\cdot\sigma^{-1}(\theta)\frac{X_\theta-Y_\theta}{\zeta(\theta)}\d\theta\right)(r)\right|^2\d r\Bigg]\cr
&=&\E_t\exp\Bigg[-\frac{1}{p-1}\int_0^t\left\langle K_H^{-1}\left(\int_0^\cdot\sigma^{-1}(\theta)\frac{X_\theta-Y_\theta}{\zeta(\theta)}\d\theta\right)(r),\d\tilde{W}_r\right\rangle\cr
&&-\frac{1}{2}\frac{1}{(p-1)^2}\int_0^t\left| K_H^{-1}\left(\int_0^\cdot\sigma^{-1}(\theta)\frac{X_\theta-Y_\theta}{\zeta(\theta)}\d\theta\right)(r)\right|^2\d r\cr
&&+\frac{1}{2}\frac{p}{(p-1)^2}\int_0^t\left| K_H^{-1}\left(\int_0^\cdot\sigma^{-1}(\theta)\frac{X_\theta-Y_\theta}{\zeta(\theta)}\d\theta\right)(r)\right|^2\d r\Bigg]\cr
&\leq&\exp\Bigg[\frac{p}{(p-1)^2}(C+C'T+C''T^{2\alpha_0})\frac{T^{2(1-H)}}{(1-\e^{-\frac{2}{3}K T})^3}|x-y|^2\Bigg].
\eeqnn
This, together with the Fatou lemma and \eqref{3.17}, leads to the desired Harnack inequality.
\fin

\bremark\label{R3.1}
By a Lamperti transform and Theorem \ref{T3.1}, we can derive Harnack type inequalities for one-dimensional SDE by multiplicative noise with $H>1/2$:
\beqnn
\d X_t=b(t,X_t)\d t+\sigma(t,X_t)\d B_t^H,\ X_0=x\in\R, \ t\in[0,T].
\eeqnn
\eremark

As a direct application of the Harnack type inequalities derived above, by \cite[Proposition 4.1]{Prato&Rockner&Wang09a} we get the strong Feller property on $P_T$.
\bcorollary\label{C3.1}
Assume (H1).
Then $P_T$ is strong Feller, i.e. the relation
\beqnn
\lim\limits_{|y-x|\rightarrow 0}P_T f(y)=P_T f(x).
\eeqnn
holds for each $f\in\mathscr{B}_b(\R^d)$ and $x\in\R^d$.
\ecorollary

To present some more applications of Theorem \ref{T3.1}, let us introduce some notations and another assumption.

Observe that the solution $X$ of equation \eqref{2.1} is not a Markov process.
As a consequence, $(P_T)_{T\geq0}$ does not consist of a Markov semigroup.
Thus, we construct the following semigroup in discrete time, i.e. for any Borel set $A$ in $\R^d$,
\beqnn P_T(x,A):=P_TI_A(x),\ \
P_T^n(x,A):=\int_{\mathbb{R}^d}P_T^{n-1}(x,\d y)P_T(y,A), \ n\geq 2.
\eeqnn
In general, $(P_T^nf)(x)=\int_{\mathbb{R}^d} f(y)P_T^n(x,\d y),\ x\in\R^d,\ f\in\mathscr{B}_b(\R^d)$.\\
Next, we shall focus on the existence and uniqueness of invariant probability measure for the discrete semigroup $(P_T^n)_{n\geq1}$,
and if so, discuss its entropy-cost inequality.
To this end, we assume moreover (H2)
\begin{itemize}
\item[] $\langle x,b(t,x)\rangle\leq L|x|^2, \ \ \forall t\in[0,T], x\in\R^d$,
\end{itemize}
where $L\in\R$ satisfies $Me^{2LT}<1,\ M$ is a positive constant given in Lemma \ref{L3.2} below.

\btheorem\label{T3.2}
Assume (H1) and (H2).
Then the semigroup $(P_T^n)_{n\geq1}$ has a unique invariant measure $\mu$.
\etheorem

In order to verify this theorem, a preliminary estimate is necessary.

\blemma\label{L3.2}
Assume (H1) and (H2).
Then there exists a positive constant $M$ (independence of $L$) such that
\beqnn
\E|X_T^x|^2\leq M\e^{2LT}(1+|x|^2).
\eeqnn
\elemma

\emph{Proof.}
According to the change-of-variables formula \cite[Theorem 4.3.1]{Zahle98a} and (H2), we get
\beqlb\label{3.18}
|X_T^x|^2\ar=\ar|x|^2+2\int_0^T\langle X_t^x,b(t,X_t^x)\rangle\d t+2\int_0^T\langle \sigma^*(t)X_t^x,\d B_t^H\rangle\cr
         \ar\leq\ar|x|^2+2L\int_0^T|X_t^x|^2\d t+2\int_0^T\langle\sigma^*(t)X_t^x,\d B_t^H\rangle.
\eeqlb
Next, we are to estimate the term $\int_0^T\langle\sigma^*(t)X_t^x,\d B_t^H\rangle$.\\
Due to the fractional integration by parts formula \eqref{2.0}, the above Riemann-Stieltjes integral can then be expressed as
\beqlb\label{3.19}
\int_0^T\langle\sigma^*(t)X_t^x,\d B_t^H\rangle=(-1)^\alpha\int_0^TD^\alpha_{0+}(\sigma^*(\cdot)X_\cdot^x)(r)D^{1-\alpha}_{T-}B^H_{T-}(r)\d r,
\eeqlb
where $1-H<\alpha<\alpha_0$, $D^\alpha_{0+}$ and $D^{1-\alpha}_{T-}$ are given by, respectively,
\beqlb\label{3.20}
D^\alpha_{0+}(\sigma^*(\cdot)X_\cdot^x)(r)=
\frac{1}{\Gamma(1-\alpha)}\left(\frac{\sigma^*(r)X_r^x}{r^\alpha}+\alpha\int_0^r\frac{\sigma^*(r)X_r^x-\sigma^*(s)X_s^x}{(r-s)^{1+\alpha}}\d s\right)
\eeqlb
 and
\beqlb\label{3.21}
D^{1-\alpha}_{T-}B^H_{T-}(r)=
\frac{(-1)^{1-\alpha}}{\Gamma(\alpha)}\left(\frac{B^H_r-B^H_T}{(T-r)^{1-\alpha}}+(1-\alpha)\int_r^T\frac{B^H_r-B^H_s}{(s-r)^{2-\alpha}}\d s\right).
\eeqlb
By \eqref{3.21}, we get, for $H>\lambda>1-\alpha$,
\beqlb\label{3.22}
|D^{1-\alpha}_{T-}B^H_{T-}(r)|\leq c\|B^H\|_\lambda(T-r)^{\lambda+\alpha-1},
\eeqlb
where and in what follows, $c$ denotes a generic constant. \\
On the other hand, by \eqref{3.20} and noting the fact that $\sigma$ is also H\"{o}lder continuous of order $\alpha_0$, we arrive at
\beqlb\label{3.23}
|D^\alpha_{0+}(\sigma^*(\cdot)X_\cdot^x)(r)|
\leq c\left(\|X^x\|_\infty r^{-\alpha}+\|X^x\|_\infty r^{\alpha_0-\alpha}+\int_0^r\frac{|X_r^x-X_s^x|}{(r-s)^{1+\alpha}}\d s\right).
\eeqlb
Observe that by the fractional integration by parts formula \eqref{2.0}, (i) of (H1) and the Gronwall lemma, we conclude that
\begin{itemize}
\item[] $\|X^x\|_\infty\leq c(1+|x|+\|B^H\|_\lambda)$,
\item[] $|X_r^x-X_s^x|\leq c\left[(1+|x|)|r-s|+\|B^H\|_\lambda(|r-s|+|r-s|^\lambda+|r-s|^{\lambda+\alpha_0})\right]$.
\end{itemize}
Substituting the two previous estimates into \eqref{3.23} yields
\beqlb\label{3.24}
&&|D^\alpha_{0+}(\sigma^*(\cdot)X_\cdot^x)(r)|\cr
&\leq&c\Big[(1+|x|)\left(r^{-\alpha}+r^{1-\alpha}+r^{\alpha_0-\alpha}\right)\cr
&&+\|B^H\|_\lambda\left(r^{-\alpha}+r^{1-\alpha}+r^{\alpha_0-\alpha}+r^{\lambda-\alpha}+r^{\lambda+\alpha_0-\alpha}\right)\Big].
\eeqlb
Combining \eqref{3.19}, \eqref{3.22} with \eqref{3.24}, we obtain
 \beqnn
\left|\int_0^T\langle\sigma^*(t)X_t^x,\d B_t^H\rangle\right|\leq c(1+|x|^2+\|B^H\|_\lambda^2).
 \eeqnn
This, together with \eqref{3.18} and the Gronwall lemma, confirms the assertion. \fin

Now, we proceed with the proof of Theorem \ref{T3.2}.

\emph{Proof of Theorem \ref{T3.2}.} Existence:
We will make use of Krylov-Bogoliubov's method.\\
Let $x_0\in\mathbb{R}^d$ and define
\beqnn
\mu_n:=\frac{\sum\limits_{k=1}^n\delta_{x_0}P_T^k}{n},\ \ n\geq1,
\eeqnn
i.e. for each $f\in\mathscr{B}_b(\mathbb{R}^d),\ \mu_n(f)=\frac{\sum_{k=1}^nP_T^kf(x_0)}{n}$.\\
Next, we will prove the tightness of $\{\mu_n\}_{n\geq1}$.\\
Firstly, based on induction argument we shall show that $\left\{\int_{\mathbb{R}^d}|x|^2P_T^n(x_0,\d x)\right\}_{n\geq1}$ is bounded. \\
When $n=1$, it follows directly from Lemma \ref{L3.2} that
\beqnn
\int_{\R^d}|x|^2P_T(x_0,\d x)\leq M\e^{2LT}(1+|x_0|^2)=:a(1+|x_0|^2).
\eeqnn
Suppose that
\beqnn
\int_{\R^d}|x|^2P_T^{n-1}(x_0,\d x)\leq a+a^2+\cdot\cdot\cdot+a^{n-1}+a^{n-1}|x_0|^2
\eeqnn
holds, then by Lemma \ref{L3.2} again we obtain
\beqnn
\int_{\R^d}|x|^2P_T^n(x_0,\d x)\ar=\ar\int_{\R^d}\int_{\R^d}|x|^2P_T^{n-1}(x_0,\d y)P_T(y,\d x)
=\int_{\R^d}\E|X_T^y|^2P_T^{n-1}(x_0,\d y)\cr
\ar\leq\ar a+a\int_{\R^d}|y|^2P_T^{n-1}(x_0,\d y)\cr
\ar\leq\ar a+a\left(a+a^2+\cdot\cdot\cdot+a^{n-1}+a^{n-1}|x_0|^2\right)\cr
\ar=\ar a+a^2+\cdot\cdot\cdot+a^n+a^n|x_0|^2.
\eeqnn
Hence, we have, for any $n\geq1$,
\beqnn
\int_{\R^d}|x|^2P_T^n(x_0,\d x)\leq\frac{a}{1-a}+|x_0|^2. \eeqnn
Consequently, there holds
\beqnn
\int_{\R^d}|x|^2\mu_n(\d x)=\frac{\sum\limits_{k=1}^n\int_{\R^d}|x|^2P_T^k(x_0,\d x)}{n}\leq\frac{a}{1-a}+|x_0|^2.
\eeqnn
Using the Chebyshev inequality, we have
\beqnn
\sup\limits_n\mu_n(|\cdot|^2>r)\leq\frac{1}{r}\left(\frac{a}{1-a}+|x_0|^2\right)\rightarrow0, \ r\rightarrow\infty,
\eeqnn
which shows the tightness of $\{\mu_n\}_{n\geq1}$.\\
So, from the Prohorov theorem, there exists a probability $\mu$ and a subsequence ${\mu_{n_k}}$ such that
$\mu_{n_k}\rightarrow\mu$ weakly as $n_k\rightarrow\infty$.
To simplify  notation, we will denote $\mu_n\rightarrow\mu$ weakly as $n\rightarrow\infty$.
Now, we will prove that $\mu$ is a invariant probability measure for $(P_T^n)_{n\geq1}$.\\
Denote $C_b(\mathbb{R}^d)$ by the set of all bounded continuous functions on $\mathbb{R}^d$. \\
For any $f\in C_b(\mathbb{R}^d)$, it follows from Corollary \ref{C3.1} that $P_T^nf\in C_b(\R^d),\ \forall n\geq1$.
Furthermore, we conclude that, for all $l\in\mathbb{N}$,
\beqnn
\mu(P_T^lf)\ar=\ar\lim\limits_{n\rightarrow\infty}\mu_n(P_T^lf)\cr
\ar=\ar\lim\limits_{n\rightarrow\infty}\frac{\sum\limits_{k=1}^nP_T^{k+l}f(x_0)}{n}\cr
\ar=\ar\lim\limits_{n\rightarrow\infty}\frac{\sum\limits_{m=1}^nP_T^{m}f(x_0)}{n}
+\lim\limits_{n\rightarrow\infty}\frac{\sum\limits_{m=n+1}^{n+l}P_T^{m}f(x_0)}{n}
-\lim\limits_{n\rightarrow\infty}\frac{\sum\limits_{m=1}^lP_T^{m}f(x_0)}{n}\cr
\ar=\ar\mu(f), \eeqnn
i.e. $\mu$ is invariant for $(P_T^n)_{n\geq1}$ .

Uniqueness: By Theorem \ref{T3.1} and \cite[Theorem 1.4.1]{Wang13a}, the proof is completed. \fin

To conclude this section, we present below the entropy-cost inequality for $\mu$.

\bcorollary\label{C3.2}
Assume (H1) and (H2).
Then for the above invariant measure $\mu$, the entropy-cost inequality
\beqnn
\mu(P_T^*f\log P_T^*f)\leq(C+C'T+C''T^{2\alpha_0})\frac{T^{2(1-H)}}{(1-\e^{-\frac{2}{3}K T})^3}W_2^2(\mu,f\mu) \eeqnn
holds for non-negative $f\in\mathscr{B}_b(\R^d)$ with $\mu(f)=1$,
where $P_T^*$ is the adjoint operator of $P_T$ in $L^2(\mu)$
and $W_2$ is the $L^2$-Wasserstein distance induced by the Euclidian metric,
i.e. for any two probability measures $\mu_1,\mu_2$ on $\R^d$,
\beqnn
W_2^2(\mu_1,\mu_2):=\inf\limits_{\pi\in\mathscr{C}(\mu_1,\mu_2)}\int_{\R^d\times\R^d}|x-y|^2\pi(\d x,\d y), \eeqnn
where $\mathscr{C}(\mu_1,\mu_2)$ is the set of all couplings of $\mu_1$ and $\mu_2$.
\ecorollary

\emph{Proof.}
Applying Theorem \ref{T3.1} to $P_T^*f$ in place of $f$, we have
\beqlb\label{3.25}
P_T(\log P_T^*f)(y)\leq\log P_T(P_T^*f)(x)+(C+C'T+C''T^{2\alpha_0})\frac{T^{2(1-H)}}{(1-\e^{-\frac{2}{3}K T})^3}|x-y|^2,\ x,y\in\R^d.
\eeqlb
Integrating both sides of \eqref{3.25} with respect to $\pi\in\mathscr{C}(\mu,f\mu)$ yields
\beqlb\label{3.26}
\mu(P_T^*f\log P_T^*f)&\leq&\mu( \log P_T(P_T^*f))\cr
&&+(C+C'T+C''T^{2\alpha_0})\frac{T^{2(1-H)}}{(1-\e^{-\frac{2}{3}K T})^3}\int_{\R^d\times\R^d}|x-y|^2\pi(\d x,\d y).
\eeqlb
Observe that from the Jensen inequality and invariance of $\mu$, we have
\beqnn
\mu( \log P_T(P_T^*f))\leq\log\mu(P_T(P_T^*f))=\log\mu(P_T^*f)=\log\mu(f P_T1)=\log\mu(f)=0.
\eeqnn
Therefore, \eqref{3.26} becomes
\beqnn
\mu(P_T^*f\log P_T^*f)\leq(C+C'T+C''T^{2\alpha_0})\frac{T^{2(1-H)}}{(1-\e^{-\frac{2}{3}K T})^3}\inf\limits_{\pi\in\mathscr{C}(\mu,f\mu)}\int_{\R^d\times\R^d}|x-y|^2\pi(\d x,\d y).
\eeqnn
This completes the proof. \fin

\textbf{Acknowledgement}
The author would like to thank Professor Feng-Yu Wang for his encouragement and comments that have led to
improvements of the manuscript.


\begin{thebibliography}{99}

\bibitem{Aida&Zhang01a} S. Aida and H. Kawabi, {\it Short time asymptotics of a certain infinite dimensional diffusion process},
  Stochastics Analysis and Related Topics 48(2001), 77--124.

\bibitem{Aida&Zhang02a}  S. Aida and T. Zhang, {\it On the small time asymptotics of diffusion processes on path groups}, Potential Anal. 16(2002), 67--78.

\bibitem{Alos&Mazet&Nualart01a} E. Al\`{o}s, O. Mazet and D. Nualart,  {\it Stochastic calculus with respect to Gaussian processes}, Ann. Probab.
  29(2001), 766--801.

\bibitem{Biagini&Hu08a} F. Biagini, Y. Hu, B. $\emptyset$ksendal and  T. Zhang, {\it Stochastic Calculus for Fractional Brownian Motion and Applications}, Springer-Verlag, London, 2008.

\bibitem{Bobkov&Gentil&Ledoux01a} S. Bobkov, I. Gentil and M. Ledoux, {\it Hypercontractivity of Hamilton-Jacobi equations}, J. Math. Pures Appl. 80(2001), 669--696.

\bibitem{Coutin&Qian02a} L. Coutin and Z. Qian, {\it Stochastic analysis, rough path analysis and fractional Brownian motions},
   Probab. Theory Related Fields 122(2002), 108--140.

\bibitem{Decreusefond&Ustunel98a} L. Decreusefond and A. S. \"{U}st\"{u}nel, {\it Stochastic analysis of the fractional Brownian motion}, Potential Anal.
    10(1998), 177--214.

\bibitem{Fan12a} X. L. Fan, {\it Derivative formula, integration by parts formula and applications for SDEs driven by fractional Brownian motion},
   arXiv:1206.0961.

\bibitem{Fan13a} X. L. Fan, {\it Harnack inequality and derivative formula for SDE driven by fractional Brownian motion},
  Science in China-Mathematics 561(2013), 515--524.

\bibitem{Fan13b} X. L. Fan, {\it Derivative formulae and Harnack inequalities for SDEs with fractional noises},  arXiv:1308.5309.

\bibitem{Gong&Wang01a} F. Gong and F. Y. Wang, {\it Heat kernel estimates with application to compactness of manifolds}, Q. J. Math. 52(2001), 171--180.

\bibitem{Guillin&Wang11a} A. Guillin and F. Y. Wang, {\it Degenerate Fokker-Planck equations : Bismut formula, gradient estimate and Harnack inequality}, J. Differential Equations 253(2012), 20--40.

\bibitem{Kawabi05a} H. Kawabi, {\it The parabolic Harnack inequality for the time dependent Ginzburg-Landau type SPDE and its application}, Potential Anal. 22(2005), 61--84.

\bibitem{Liu09a} W. Liu,  {\it Harnack inequality and applications for stochastic evolution equations with monotone drifts}, J. Evol. Equ. 9(2009), 747--770.

\bibitem{Lyons98a} T. Lyons, {\it Differential equations driven by rough signals}, Rev. Mat. Iberoamericana 14(1998), 215--310.

\bibitem{Nikiforov&Uvarov88} A. F. Nikiforov and V. B. Uvarov, {\it Special Functions of Mathematical Physics}, Birkh\"{a}user, Boston, 1988.

\bibitem{Nualart06a} D. Nualart, {\it The Malliavin Calculus and Related Topics, Second edition}, Springer-Verlag, Berlin, 2006.

\bibitem{Nualart&Ouknine02b} D. Nualart and Y. Ouknine,  {\it Regularization of differential equations by fractional noise},  Stochastic Process. Appl.
    102(2002), 103--116.
\bibitem{Nualart&Rascanu02a} D. Nualart and A. R\u{a}\c{s}canu, {\it Differential equations driven by fractional Brownian motion}, Collect. Math. 53(2002), 55--81.
\bibitem{Ouyang09a} S. X. Ouyang,  {\it Harnack inequalities and applications for stochastic equations}, Ph.D. Thesis, Bielefeld University, 2009.

\bibitem{Ouyang&Rockner&Wang12a} S. X. Ouyang,  M. R\"{o}ckner and F. Y. Wang,  {\it Harnack inequalities and applications for Ornstein-Uhlenbeck semigroups with jump}, Potential Anal. 36(2012), 301--315.

\bibitem{Prato&Rockner&Wang09a} G. Da Prato, M. R\"{o}ckner and F.Y. Wang, {\it Singular stochastic equations on Hilberts space: Harnack inequalities for their transition semigroups}, J. Funct. Anal. 257(2009), 992--1017.

\bibitem{Rockner&Wang03a} M. R\"{o}ckner and F. Y. Wang, {\it Supercontractivity and ultracontractivity for (non-symmetric) diffusion semigroups on manifolds}, Forum Math. 15(2003), 893--921.

\bibitem{Rockner&Wang03b} M. R\"{o}ckner and F. Y. Wang, {\it Harnack and functional inequalities for generalized Mehler semigroups}, J. Funct. Anal. 203(2003), 237--261.

\bibitem{Rockner&Wang10a} M. R\"{o}ckner and F. Y. Wang, {\it Log-Harnack inequality for stochastic differential equations in
Hilbert spaces and its consequences}, Infin. Dimens. Anal. Quantum Probab. Relat. Top. 13(2010), 27--37.

\bibitem{Samko&Kilbas&Marichev} S. G. Samko, A. A. Kilbas and O. I. Marichev,  {\it Fractional Integrals and Derivatives, Theory and Applications}, Gordon
    and Breach Science Publishers, Yvendon, 1993.

\bibitem{Wang97a} F. Y. Wang, {\it Logarithmic Sobolev inequalities on noncompact Riemannian manifolds},
 Probab. Theory Related Fields 109(1997), 417--424.

\bibitem{Wang99a}  F. Y. Wang, {\it Harnack inequalities for log-Sobolev functions and estimates of log-Sobolev constants}, Ann. Probab. 27(1999), 653--663.

\bibitem{Wang07a}  F. Y. Wang, {\it Harnack inequality and applications for stochastic generalized
                  porous media equations}, Ann. Probab. 35(2007), 1333--1350.

\bibitem{Wang10a} F. Y. Wang, {\it Harnack inequalities on manifolds with boundary and applications}, J. Math. Pures Appl. 94(2010), 304--321.

\bibitem{Wang11b} F.Y. Wang,  {\it Harnack inequality for SDE with multiplicative noise and extension to Neumann semigroup on nonconvex manifolds}, Ann. Probab. 39(2011), 1449--1467.

\bibitem{Wang13a} F. Y. Wang, {\it Harnack Inequalities for Stochastic Partial Differential Equations}, Springer, 2013.

\bibitem{Wang&Yuan11a}  F. Y. Wang and C. Yuan, {\it Harnack inequalities for functional SDEs with multiplicative noise
                    and applications}, Stochastic Process. Appl. 121(2011), 2692--2710.

\bibitem{Young36a} L. C. Young,  {\it An inequality of the H\"{o}lder type connected with Stieltjes integration}, Acta Math. 67(1936), 251--282.


\bibitem{Zahle98a} M. Z\"{a}hle, {\it Integration with respect to fractal functions and stochastic calculus I}, Probab. Theory Related Fields 111(1998),
    333--374.

\bibitem{Zhang13a}  X. C. Zhang, {\it Derivative formulas and gradient estimates for SDEs driven by $\alpha$-stable processes},
Stochastic Process. Appl. 123(2013), 1213--1228.




\end{thebibliography}
\end{document}